\RequirePackage{ifpdf}
\ifpdf % We are running pdfTeX in pdf mode
\documentclass[pdftex]{sigma}
\else
\documentclass{sigma}
\fi

\def\frak{\mathfrak}
\def\Bbb{\mathbb}
\def\Cal{\mathcal}

\newtheorem*{prop*}{Proposition}

\newtheorem*{thm*}{Theorem}

\newtheorem*{lem*}{Lemma}

\newtheorem*{kor*}{Corollary}

\newcommand{\ad}{\operatorname{ad}}
\newcommand{\Ad}{\operatorname{Ad}}

\newcommand{\id}{\operatorname{id}}

\renewcommand{\ker}{\operatorname{ker}}
\newcommand{\im}{\operatorname{im}}

\newcommand{\fg}{{\mathfrak g}}

\newcommand{\x}{\times}
\renewcommand{\o}{\circ}

\let\ccdot\cdot
\def\cdot{\hbox to 2.5pt{\hss$\ccdot$\hss}}

\newcommand{\al}{\alpha}

\newcommand{\la}{\lambda}

\newcommand{\om}{\omega}
\renewcommand{\phi}{\varphi}
\newcommand{\ph}{\varphi}

\newcommand{\si}{\sigma}

\newcommand{\Ga}{\Gamma}
\newcommand{\La}{\Lambda}
\newcommand{\Om}{\Omega}
\newcommand{\Ph}{\Phi}

\begin{document}

\allowdisplaybreaks
	
\renewcommand{\PaperNumber}{111}

\FirstPageHeading

\renewcommand{\thefootnote}{$\star$}

\ShortArticleName{Curved Casimir Operators and the BGG Machinery}

\ArticleName{Curved Casimir Operators and the BGG Machinery\footnote{This paper is a
contribution to the Proceedings of the 2007 Midwest
Geometry Conference in honor of Thomas~P.\ Branson. The full collection is available at
\href{http://www.emis.de/journals/SIGMA/MGC2007.html}{http://www.emis.de/journals/SIGMA/MGC2007.html}}}

\Author{Andreas \v CAP~$^{\dag \ddag}$ and Vladim\'{i}r SOU\v CEK~$^\S$}

 \AuthorNameForHeading{A. \v Cap and V. Sou\v cek}

\Address{$^{\dag}$~Fakult\"at f\"ur Mathematik, Universit\"at Wien,
  Nordbergstr.~15, A-1090 Wien, Austria}

\Address{$^{\ddag}$~International Erwin
  Schr\"odinger Institute for Mathematical Physics,\\
$\phantom{^{\ddag}}$~Boltzmanngasse 9,
  A-1090 Wien, Austria}
\EmailD{\href{mailto:Andreas.Cap@esi.ac.at}{Andreas.Cap@esi.ac.at}}

\Address{$^\S$~Mathematical Institute, Charles University,
  Sokolovsk\'a 83, Praha, Czech Republic}
\EmailD{\href{mailto:soucek@karlin.mff.cuni.cz}{soucek@karlin.mff.cuni.cz}}

\ArticleDates{Received August 24, 2007, in f\/inal form November
16, 2007; Published online November 22, 2007}

\Abstract{We prove that the Casimir operator acting on sections of a
  homogeneous vector bundle over a generalized f\/lag manifold naturally
  extends to an invariant dif\/ferential opera\-tor on arbitrary parabolic
  geometries. We study some properties of the resulting invariant
  operators and compute their action on various special types of
  natural bundles. As a f\/irst application, we give a very general
  construction of splitting operators for parabolic geometries. Then
  we discuss the curved Casimir operators on dif\/ferential forms with
  values in a~tractor bundle, which nicely relates to the machinery of
  BGG sequences. This also gives a nice interpretation of the
  resolution of a f\/inite dimensional representation by (spaces of
  smooth vectors in) principal series representations provided by a
  BGG sequence.}

\Keywords{induced representation; parabolic geometry; invariant
  dif\/ferential operator; Casimir operator; tractor bundle; BGG sequence}

\Classification{22E46; 53A40; 53C15; 58J70}

\rightline{\textit{Dedicated to the memory of Thomas P. Branson}}

\section{Introduction}\label{1}
Let $G$ be a semisimple Lie group and $P\subset G$ a parabolic
subgroup. Then to any representation of $P$ one may associate a
representation of $G$ by the process of parabolic induction. The
induced representations are most naturally def\/ined on spaces of
sections of homogeneous vector bundles on the generalized f\/lag manifold
$G/P$.  As on any representation of $G$, the Casimir operator
naturally acts on (the smooth vectors of) any induced representation,
and it is easy to see that it acts by a $G$-invariant dif\/ferential
operator.

Starting from examples like conformal structures and CR-structures,
it has been realized that manifolds endowed with certain geometric
structures can be viewed as ``curved analogs'' of homogeneous spaces
of the form $G/P$. More formally, they can be equivalently described
as certain Cartan geometries with homogeneous model $G/P$. These
structures are called parabolic geometries of type $(G,P)$ and they
have been intensively studied during the last years.

Any representation of $P$ gives rise to a natural vector bundle on all
parabolic geometries of type $(G,P)$, and on the homogeneous model,
this gives rise to the usual homogeneous vector bundles.
The concept of $G$-invariant dif\/ferential operators between sections
of homogeneous vector bundles generalizes to the notion of invariant
(or natural) dif\/ferential operators on parabolic geometries. Questions
on such operators are surprisingly subtle and are a major research
topic, see e.g.~\cite{Eastwood-Rice, GJMS, Graham, CSS-BGG,
  David-Tammo, Gover-Hirachi, Gover-Graham}.

The f\/irst main result of this article is that the Casimir operator on
an induced representation canonically extends to a natural operator on
arbitrary parabolic geometries of type $(G,P)$, which we call the
\textit{curved Casimir operator}. This extension is described in terms
of the so-called fundamental derivative, which is easily connected to
well developed tools for parabolic geometries such as tractor calculus.

Next, we prove that, similarly to the original Casimir operator, the
curved analog has strong naturality properties. Using ideas closely
related to tractor calculus we derive a simple formula for the curved
Casimir operator on an arbitrary natural bundle. This is used to show
that the curved Casimir always is an operator of order at most one. On
natural bundles induced by irreducible representations of $P$ (which
correspond to principal series representations) the curved Casimir
acts by a scalar, and we compute this scalar in terms of weights. We
sketch how this result can be applied to construct a large number of
natural dif\/ferential operators.

In the last section, we consider the curved Casimir operators on the
so--called tractor bundles, see \cite{TAMS}, as well as on tractor
bundle valued dif\/ferential forms. On the homogeneous model, tractor
bundles correspond to representations induced by restrictions to $P$
of representations of $G$, and they are canonically trivial. In the
curved case, however, they need not be trivial and they are important
tools, since they admit canonical linear connections. We related the
curved Casimir operator to these tractor connections and the induced
covariant exterior derivative.

A large number of invariant dif\/ferential operators on parabolic
geometries can be constructed using the machinery of
Bernstein--Gelfand--Gelfand (or BGG) sequences, which applies to
curved geometries. On the homogeneous model, one obtains resolutions
of all f\/inite dimensional irreducible representation of $G$ by
invariant dif\/ferential operators between sections of homogeneous
vector bundles.  These operators are related by a duality to the
homomorphisms of generalized Verma modules in the generalized BGG
resolution of the dual of the given f\/inite dimensional representation.
BGG sequences were f\/irst constructed in \cite{CSS-BGG}, and the
construction was improved in \cite{David-Tammo}. Our description of
the curved Casimir operator sheds new light on this construction. In
particular, we prove that the resolution on the homogeneous model is
obtained by restricting the twisted de-Rham resolution to eigenspaces
of the Casimir operator.

The interaction between results for the homogeneous model (which often
are deduced using representation theory) and results on general curved
geometries is one of the basic features of the theory of parabolic
geometries. This has led to interesting interactions between people
working in representation theory and people working on parabolic
geometries, but the communication between the two groups is not always
easy. One of the few people active in both communities was Tom
Branson, who did a lot to facilitate mutual communication. He died
completely unexpectedly in 2006, and the 2007 Midwest Geometry
Conference was held in his honor. This article is our contribution to
the proceedings of the MGC 2007 conference. Trying to make the article
accessible for people from both communities, we have included a bit of
review material as well as comparison between the languages and
notations which are standardly used in the two f\/ields. We hope that
this article will be useful for intensifying the interactions between
the two communities.

\section{A curved analog of the Casimir operator}\label{2}

\subsection[Generalized flag manifolds and induced representations]{Generalized f\/lag manifolds and induced representations}\label{2.1}

A generalized f\/lag manifold is the quotient $G/P$ of a semisimple Lie
group $G$ by a parabolic subgroup $P\subset G$. These homogeneous
spaces are always compact, and in the case that $G$ is complex, they
are the only compact homogeneous spaces of $G$. In representation
theory, these spaces play an important role via parabolic induction.
Suppose that $W$ is a f\/inite dimensional representation of $P$. Then
one forms the homogeneous vector bundle $G\x_PW=(G\x W)/P$, where $P$
acts from the right on $G\x W$ via $(g,w)\cdot h:=(gh,h^{-1}\cdot w)$.
The $G$-action on $G/P$ naturally lifts to an action on $G\x_PW$ by
vector bundle homomorphisms. Hence the space $\Ga(G\x_PW)$ of smooth
sections carries a natural representation of $G$, called the
\textit{representation induced by~$W$}, def\/ined by $(g\cdot s)(\tilde
gP)=g\cdot (s(g^{-1}\tilde gP))$. Using compactness of $G/P$, this
representation can be extended to various completions of
$\Ga(G\x_PW)$, and this construction is one of the basic sources of
representations of $G$. If one starts with an irreducible
representation of $P$, the resulting representation of $G$ is called a
\textit{principal series} representation.

There is an equivalent description of induced representations which is
commonly used in representation theory. The space $\Ga(G\x_PW)$ can be
naturally identif\/ied with
\[
C^\infty(G,W)^P=\{f\in C^\infty(G,W):f(gh)=h^{-1}\cdot
f(g)\quad\forall\, g\in G,h\in P\}.
\]
(This depends only on the fact that $G\x_PW$ is associated to the
principal bundle $G\to G/P$.) The correspondence between
$s\in\Ga(G\x_PW)$ and $f\in C^\infty(G,W)^P$ is characterized by the
fact that $s(gP)$ is the equivalence class of $(g,f(g))$. In this
picture, the natural left action of $G$ is given by $(g\cdot f)(\tilde
g)=f(g^{-1}\tilde g)$.

In some cases, the $G$-action on $G/P$ admits a direct geometric
interpretation. For example, if $G=O(n+1,1)$ and $P\subset G$ is the
stabilizer of an isotropic line, then $G/P$ is the $n$-sphere~$S^n$,
and for $n\geq 3$, the action on $G/P$ identif\/ies $G$ with the group
of conformal isometries of the round metric on $S^n$. These are
exactly the dif\/feomorphisms of $S^n$, which only rescale the round
metric by a positive smooth function. For $G=SU(n+1,1)$ and $P\subset
G$ the stabilizer of a~complex isotropic line, one obtains
$G/P=S^{2n+1}$, viewed as a hypersurface in $\Bbb C^{n+1}$. The action
on~$G/P$ identif\/ies $G$ with the group of all CR-dif\/feomorphisms of
$S^{2n+1}$. These are exactly the dif\/feomorphisms induced from
biholomorphisms of the unit ball.

\subsection{Parabolic geometries}\label{2.2}

It is an idea due to Elie Cartan to describe certain geometric
structures as ``curved analogs'' of homogeneous spaces, which
represent geometries in the sense of Klein's Erlangen program. An
exposition of this approach and several examples can be found in the
book \cite{Sharpe} and in \cite{Srni05}, which is specif\/ically
directed towards the parabolic case. The starting point of this
approach is that the left actions of the elements of $G$ on a
homogeneous space $G/H$ can be characterized as those automorphisms of
the $H$-principal bundle $G\to G/H$ which preserve the left
Maurer--Cartan form on $G$. A curved analog, called a Cartan geometry
of type $(G,H)$, is then given by a principal $H$-bundle $\Cal G$ on
a manifold $M$ such that $\dim(M)=\dim(G/H)$ and a Cartan connection
$\om\in\Om^1(\Cal G,\frak g)$, where $\frak g$ is the Lie algebra of
$G$. The Cartan connection is required to def\/ine an isomorphism $T\Cal
G\to\Cal G\x\frak g$, to be equivariant with respect to the principal
right action of $H$, and to reproduce the generators of fundamental
vector f\/ields.

Cartan himself proved that conformal structures in dimension $\geq 3$
and CR structures in dimension $3$ can be equivalently described as
Cartan geometries corresponding to the homogeneous spaces of
$O(n+1,1)$ respectively $SU(2,1)$ mentioned above. For CR structures
in arbitrary dimensions this is a famous result of Chern and Moser,
see \cite{CM}, with an earlier proof due to Tanaka. The fact that the
Cartan geometries provide an equivalent description of the underlying
structure is a non-trivial theorem in these cases. In pioneering work
culminating in \cite{Tanaka79}, N.~Tanaka showed that for every
semisimple Lie group $G$ and parabolic subgroup $P\subset G$, Cartan
geometries of type $(G,P)$ (which satisfy the additional conditions of
regularity and normality) are equivalent encodings of certain
underlying structures.  The description of these underlying structures
has been later simplif\/ied considerably, see \cite{Srni05}.  The
geometric structures obtained in this way have been named parabolic
geometries and have been intensively studied during the last years.
Apart from conformal and CR structures, also projective and
quaternionic structures, path geometries, quaternionic contact
structures (see \cite{Biquard}) and various types of generic
distributions are examples of parabolic geometries.

One basic principle in the study of parabolic geometries (and also of
general Cartan geometries) is the passage from $G/P$, which in this
context is referred to as the \textit{homogeneous model}, and general
(curved) Cartan geometries of type $(G,P)$.  Sometimes, this passage
is easy.  For example, homogeneous vector bundles on $G/P$ are in
bijective correspondence with representations of $P$. Given a
representation $W$ of $P$ and a Cartan geometry $(p:\Cal G\to M,\om)$
of type $(G,P)$ one can form the associated bundle $WM:=\Cal G\x_PW$.
The resulting bundles are called \textit{natural bundles} on Cartan
geometries of type $(G,P)$.

The Cartan connection $\om$ can be used to identify such bundles with
other types of natural bundles. The basic example is provided by the
representation $\frak g/\frak p$. It is well known that the
homogeneous bundle $G\x_P(\frak g/\frak p)$ is isomorphic to the
tangent bundle $T(G/P)$. Likewise, via the Cartan connection $\om$,
one constructs an isomorphism $\Cal G\x_P(\frak g/\frak p)\cong TM$.
As in the homogeneous case, sections of natural bundles can be
described via equivariant functions, i.e.~$\Ga(\Cal G\x_PW)\cong
C^\infty(\Cal G,W)^P$, with the principal right action of $P$ on $\Cal
G$ replacing the right translations by elements of $P$ on $G$.

A $G$-equivariant vector bundle map between homogeneous vector
bundles is given by a $P$-equivariant linear map between the inducing
representations and hence extends to arbitrary Cartan geometries as
well. This means, that $G$-equivariant maps between induced
representations, which are tensorial (i.e.~linear over
$C^\infty(G/P,\Bbb R)$) canonically extend to tensorial operators on
all Cartan geometries. For more general intertwining operators between
induced representations, the situation is much more complicated. In
particular, the case of intertwining dif\/ferential ope\-rators is related
to questions on invariant dif\/ferential operators on parabolic
geometries, which is a very active area of research.

\subsection{The fundamental derivative}\label{2.3}
Finding invariant dif\/ferential operators is particularly dif\/f\/icult for
bundles induced by irreducible representations of $P$. Allowing more
general bundles, there is a basic set of invariant operators. The
construction is easy to understand on the homogeneous model. Suppose
that $W$ is any representation of $P$ and consider sections of the
induced bundle $G\x_PW$ as equivariant functions $f:G\to W$. For a
vector f\/ield $\xi\in\frak X(G)$, the derivative $\xi\cdot f$ will not
be equivariant in general. This will be the case, however, if $\xi$ is
invariant under right translations by elements of $P$. Now using the
left trivialization of the tangent bundle $TG$, vector f\/ields on $G$
correspond to smooth functions $X:G\to\frak g$ via
$\xi(g)=L_{X(g)}$. Such a vector f\/ield is invariant under right
translations by elements of $P$ if and only if
$X(gh)=\Ad(h^{-1})(X(g))$ for all $g\in G$ and $h\in P$. Consequently,
$P$-invariant vector f\/ields on $G$ are in bijective correspondence
with sections of the associated bundle $G\x_P\frak g$.

This carries over to the curved case. Given a parabolic geometry
$(p:\Cal G\to M,\om)$, one def\/ines the \textit{adjoint tractor bundle}
$\Cal AM:=\Cal G\x_P\frak g$. Via the Cartan connection, smooth
sections of $\Cal AM$ are in bijective correspondence with
$P$-invariant vector f\/ields on $\Cal G$. For the natural bundle
$WM=\Cal G\x_PW$ corresponding to the $P$-representation $W$, one has
$\Ga(WM)\cong C^\infty(\Cal G,W)^P$. Dif\/ferentiating a $P$-equivariant
function with respect to a $P$-invariant vector f\/ield, the result is
$P$-equivariant again. Hence there is a bilinear natural dif\/ferential
operator
\begin{gather*}
D:\Ga(\Cal AM)\x\Ga(\Cal WM)\to\Ga(\Cal WM).
\end{gather*}
To emphasize the analogy to a covariant derivative, this operator is
usually denoted as $(s,\si)\mapsto D_s\si$. The analogy is strengthened
by the fact that $D$ is linear over $C^\infty(M,\Bbb R)$ in the $\Cal
AM$-slot and satisf\/ies a Leibniz rule in the other slot. This
operator is called the \textit{fundamental derivative} and was
introduced (under the name ``fundamental $D$-operator'') in
\cite{TAMS}. There it was also proved that this operator is compatible
with any natural bundle map. Note that we can also view $D$ as an
operator mapping $\Ga(WM)$ to $\Ga(\Cal A^*M\otimes WM)$, and in this
form the operator can be iterated without problems. Hence for a
section $\si\in\Ga(WM)$, we can form $D^k\si\in\Ga(\otimes^k\Cal
A^*M\otimes WM)$.

We have to note one more important property of the fundamental
derivative. The Lie subalgebra $\frak p\subset\fg$ is a
$P$-submodule. Hence for any parabolic geometry $(p:\Cal G\to M,\om)$
of type $(G,P)$, we obtain a subbundle $\Cal G\x_P\frak p\subset\Cal
G\x_P\frak g=\Cal AM$. For reasons that will become clear later on, we
denote this subbundle by $\Cal A^0M$. For any f\/inite dimensional
representation $W$ of $P$, we have the corresponding inf\/initesimal
representation of $\frak p$ on $W$. We can interpret this as an
action, i.e.~as a bilinear map $\frak p\x W\to W$, and this map is
$P$-equivariant. Passing to associated bundles, we obtain an induced
bilinear bundle map $\Cal A^0M\x WM\to WM$, which we write as
$(s,\si)\mapsto s\bullet\si$.

Now sections of $\Cal A^0M\subset\Cal AM$ correspond to $P$-invariant
vector f\/ields on $\Cal G$, which are mapped by $\om$ to $\frak
p\subset\fg$. By def\/inition of a Cartan connection, the value of such
a vector f\/ield in each point coincides with a fundamental vector
f\/ield, so sections of $\Cal A^0M$ correspond to $P$-invariant
vertical vector f\/ields. But dif\/ferentiating a $P$-equivariant
function with respect to a fundamental vector f\/ield, one obtains the
inf\/initesimal action of minus the generator of the fundamental f\/ield.
This implies that $D_s\si=-s\bullet\si$ holds for all $s\in\Ga(\Cal
A^0M)$ and $\si\in\Ga(WM)$.

\subsection{A curved analog of the Casimir operator}\label{2.4}
Given a (possibly inf\/inite dimensional) representation $V$ of the Lie
group $G$ and its Lie algebra $\fg$, one always has the Casimir
operator mapping the representation to itself. Since $\fg$ is
semisimple, its Killing form $B$ is non-degenerate, and hence
provides an isomorphism $\fg\cong\fg^*$ of $G$-modules. Now consider
a basis $\{\xi^\ell\}$ of $\fg$ and let $\{\xi_\ell\}$ be the dual
basis with respect to $B$. Then the Casimir operator $\Cal C:V\to V$
is def\/ined by $\Cal C(v):=\sum_\ell \xi^\ell\cdot\xi_\ell\cdot v$. It
is easy to see that $\Cal C$ is independent of the choice of the basis
$\{\xi^\ell\}$ and $G$-equivariant.

In particular, we can use this idea on the space of smooth sections of
a homogeneous bundle $G\x_PW\to G/P$. In the picture of equivariant
functions, we have seen in Section~\ref{2.1} that the $G$-action on
$C^\infty(G,W)^P$ is given by $(g\cdot f)(\tilde g):=f(g^{-1}\tilde
g)$. This immediately implies that the action can be dif\/ferentiated to
obtain an action of the Lie algebra $\fg$ of $G$. This action is
explicitly given by $X\cdot f:=-R_X\cdot f$, where $R_X\in\frak X(G)$
is the right invariant vector f\/ield generated by $X\in\fg$.

Given bases $\{\xi^\ell\}$ and $\{\xi_\ell\}$ of $\fg$ as above, we
conclude that the Casimir operator on $C^\infty(G,W)^P\!$ is given by
$\Cal C(f)=\sum_\ell R_{\xi^\ell}\cdot R_{\xi_\ell}\cdot f$. In
particular, $\Cal C$ is a dif\/ferential operator of order at most 2.
Note that from this formula, it is not evident that $\Cal C$ is
independent of the choice of basis.

We want to express the Casimir operator $\Cal C$ acting on sections of
an arbitrary homogeneous vector bundle via the fundamental derivative.
The resulting formula then def\/ines an operator on sections of the
corresponding natural vector bundle on arbitrary curved parabolic
geometries. We have noted above that the Killing form $B$ induces an
isomorphism $\fg\to\fg^*$ of $G$-modules and hence of $P$-modules.
Passing to associated bundles, we obtain a natural isomorphism $\Cal
AM\to\Cal A^*M$ between the adjoint tractor bundle and its dual on any
parabolic geometry. This in turn can be used to interpret the Killing
form as a natural tensorial linear operator which maps sections of
$\otimes^2\Cal A^*M$ to smooth functions on $M$.

\begin{proposition}\label{prop2.4}
  Let $G$ be a semisimple Lie group, $P\subset G$ a parabolic
  subgroup, let $W$ be a~representation of $P$, and consider the
  Casimir operator $\Cal C:\Ga(G\x_PW)\to\Ga(G\x_PW)$. Then~$\Cal C$
  can be written as the composition
\begin{gather*}
\Ga(W(G/P))\overset{D^2}{\longrightarrow}
\Ga(\otimes^2\Cal A^*(G/P)\otimes W(G/P))
\overset{B\otimes\id_W}{\longrightarrow} \Ga(W(G/P)).
\end{gather*}
Hence $(B\otimes\id)\o D^2$ defines an extension of $\Cal C$ to an
invariant differential operator acting on sections of the natural
bundle induced by $W$ for all parabolic geometries of type $(G,P)$.
\end{proposition}
\begin{proof}
  For any $X\in\frak g$, the right invariant vector f\/ield $R_X\in\frak
  X(G)$ is in particular invariant under right translations by
  elements of $P$. Hence it corresponds to a section $s_X\in\Ga(\Cal
  A(G/P))$, and the corresponding equivariant function $G\to\frak g$
  is given by $g\mapsto \Ad(g^{-1})X$. For a basis~$\{\xi^\ell\}$ of~$\fg$
  with dual basis~$\xi_\ell$, we thus see that $\Cal C$ can be written
  as $\Cal C(\si)=\sum_\ell D_{s_{\xi^\ell}}D_{s_{\xi_\ell}}\si$. Now
  naturality of the fundamental derivative easily implies that
  \begin{equation}
    \label{D^2}
    D^2\si(s_1,s_2)=D_{s_1}D_{s_2}\si-D_{D_{s_1}s_2}\si,
  \end{equation}
  see \cite[3.6]{TAMS}. For $X,Y\in\fg$, the function corresponding to
  $s_Y$ is $g\mapsto \Ad(g^{-1})Y$, and dif\/ferentiating this with
  respect to $R_X$, we immediately conclude that
  $D_{s_X}s_Y=s_{[Y,X]}$.

Since the value of the Casimir operator is independent of the choice of
the basis, we may as well write $\Cal C(\si)$ as
\[
  \tfrac{1}{2}\textstyle\sum_\ell( D_{s_{\xi^\ell}}D_{s_{\xi_\ell}}\si+
  D_{s_{\xi_\ell}}D_{s_{\xi^\ell}}\si),
\]
and since from above we know that $D_{s_X}s_Y$ is skew symmetric in
$X$ and $Y$, this equals
\[
\tfrac{1}{2}\textstyle\sum_\ell (D^2\si(s_{\xi^\ell},s_{\xi_\ell})+
D^2\si(s_{\xi_\ell},s_{\xi^\ell})).
\]
Invariance of the Killing form implies that
$B(\Ad(g^{-1})X,\Ad(g^{-1})Y)=B(X,Y)$, which shows that for each $g\in
P$, the elements $\{s_{\xi^\ell}(gP)\}$ and $\{s_{\xi_\ell}(gP)\}$
form dual bases for the f\/iber $\Cal A_{gP}(G/P)$ with respect to the
bundle map induced by the Killing form. Hence
\[
\textstyle\sum_\ell
D^2\si(s_{\xi^\ell},s_{\xi_\ell})=\textstyle\sum_\ell
D^2\si(s_{\xi_\ell},s_{\xi^\ell})=
(B\otimes\id)(D^2\si),
\]
and the result follows.
\end{proof}

\noindent
{\bf Def\/inition.}
We will refer to the operator $\Cal C:=(B\otimes\id)\o D^2$ on general
parabolic geometries as the \textit{curved Casimir operator}.

\section{Basic properties of the curved Casimir operators}\label{3}

\subsection{Naturality}\label{3.1}
Having a formula for the Casimir in terms of the fundamental
derivative, we can now study its basic properties not only on the
homogeneous model, but at the same time on general curved parabolic
geometries. The f\/irst consequence is that, like the fundamental
derivative, the curved Casimir operator has very strong naturality
properties.

\begin{proposition}\label{prop3.1}
  Let $W$ and $W'$ be representations of $P$, $\Ph:W\to W'$ a
  $P$-equivariant mapping, and let us also denote by $\Ph$ the
  induced natural bundle map. Denoting by $\Cal C_W$ and $\Cal C_{W'}$
  the curved Casimir operators on the natural bundles induced by $W$
  and $W'$, we get $\Cal C_{W'}(\Ph(\si))=\Ph(\Cal C_W(\si))$ for any
  parabolic geometry $(p:\Cal G\to M,\om)$ of type $(G,P)$ and any
  $\si\in\Ga(WM)$.

  In particular, the curved Casimir preserves sections of natural
  subbundles and restricts to the curved Casimir of the subbundle on
  such sections. Likewise, the induced operator on sections of a
  natural quotient bundle coincides with the curved Casimir on that
  bundle.
\end{proposition}
\begin{proof}
  By Proposition~3.1 of \cite{TAMS}, the fundamental derivative has
  the property that $D_s(\Ph(\si))=\Ph(D_s\si)$ for any $(p:\Cal G\to
  M,\om)$, $s\in\Ga(\Cal AM)$ and $\si\in\Ga(WM)$. Therefore,
  naturality follows directly from the def\/inition of $\Cal C$.

  The second part of the claim is just the specialization to the
  inclusion of a subrepresentation respectively the projection to a
  quotient.
\end{proof}

\subsection{Adapted local frames}\label{3.2}
To derive further properties of the curved Casimir, we will use a
formula in terms of appropriately chosen local dual frames for the
adjoint tractor bundle. This needs a bit more background about
parabolic subalgebras.

A basic structural result on parabolic subgroups $P\subset G$ is the
existence of a \textit{Langlands decomposition}, which is usually
written as $P=MAN$. Here $M\subset P$ is a semisimple Lie subgroup and
$A,N\subset P$ are normal subgroups which are Abelian and nilpotent,
respectively. In the literature on parabolic geometries, one usually
considers the reductive subgroup $MA\subset P$, which is commonly
denoted by $G_0$, while the nilpotent normal subgroup $N$ is usually
denoted by $P_+$. On the level of Lie algebras, one has $\frak
p_+\subset\frak p\subset\fg$ and $\frak p/\frak p_+=\fg_0$.

There is a f\/iner decomposition however, see \cite{Yamaguchi}. It turns
out that there is a grading $\fg=\fg_{-k}\oplus\dots\oplus\fg_k$ of
$\fg$ such that $\frak p=\fg_0\oplus\dots\oplus\fg_k$ and $\frak
p_+=\fg_1\oplus\dots\oplus\fg_k$. This grading is only an auxiliary
tool for parabolic geometries, since it is not $\frak p$-invariant.
The corresponding f\/iltration
\[
\fg=\fg^{-k}\supset\fg^{-k+1}\supset\dots\supset\fg^k\supset\fg^{k+1}=\{0\}
\]
def\/ined by $\fg^i:=\fg_i\oplus\dots\oplus\fg_k$ is $P$-invariant.
In this notation $\frak p=\fg^0$ and $\frak p_+=\fg^1$, so the facts
that $\frak p\subset\frak g$ is a subalgebra and that $\frak
p_+\subset\frak p$ is a nilpotent ideal immediately follow from that
fact that the f\/iltration is compatible with the Lie bracket.

For any parabolic geometry $(p:\Cal G\to M,\om)$ of type $(G,P)$, the
$P$-invariant f\/iltration of $\fg$ gives rise to a f\/iltration $\Cal
AM=\Cal A^{-k}M\supset\dots\supset\Cal A^kM$ of the adjoint tractor
bundle by smooth subbundles. In particular, $\Cal A^0M=\Cal G\x_P\frak
p$, which justif\/ies the notation introduced before. If $U\subset M$ is
an open subset over which $\Cal G$ is trivial, then all the bundles
$\Cal A^iM$ are trivial over $U$. Over such a subset, we can therefore
choose local frames of these subbundles.

From the facts that the grading and the induced f\/iltration of $\fg$
are compatible with the Lie bracket, one can easily deduce the
relation of the Killing form $B$ to these data, see again
\cite{Yamaguchi}. For $X\in\fg_i$ and $Y\in\fg_j$, $\ad_X\o\ad_Y$ maps
$\fg_\ell$ to $\fg_{i+j+\ell}$. Hence this composition is tracefree
and $B$ vanishes on $\fg_i\x\fg_j$ unless $j=-i$. Consequently, $B$ is
non-degenerate on $\fg_0$ and induces dualities of $G_0$-modules
between $\fg_i$ and $\fg_{-i}$ for all $i=1,\dots,k$. Likewise, for
$i=-k+1,\dots,k$, $\fg^i$ is the annihilator with respect to~$B$ of
$\fg^{-i+1}$, and one obtains dualities of $P$-modules between~$\fg^i$ and~$\fg/\fg^{-i+1}$.  Translated to associated bundles, this
implies that $B$ vanishes on $\Cal A^{-i+1}M\x\Cal A^iM$ and induces a
non-degenerate bilinear form on $(\Cal AM/ \Cal A^{-i+1}M)\x\Cal
A^iM$. Since $\Cal A^0M/\Cal A^1M\cong\Cal G\x_P\fg_0$, the Killing
form induces a non-degenerate bilinear form $B$ on this bundle.

\medskip

\noindent
{\bf Def\/inition.}
Consider the adjoint tractor bundle $\Cal AM$ for a parabolic geometry
of type $(G,P)$, and consider the grading and the f\/iltration of $\fg$
induced by the parabolic subalgebra $\frak p\subset\frak g$. Then an
\textit{adapted} local frame for $\Cal AM$ is a local frame of the
form
\[
\{X_i,A_r,Z^i:i=1,\dots,\dim(\frak p_+), r=1,\dots,\dim(\fg_0)\},
\]
such that
\begin{itemize}\itemsep=0pt
\item $Z^i\in\Ga(\Cal A^1M)$ for all $i$ and $A_r\in\Ga(\Cal A^0M)$
  for all $r$\newline (which implies $B(Z^i,Z^j)=0$ and $B(Z^i,A_r)=0$
  for all $i,j,r$ by the properties of $B$).
\item $B(X_i,X_j)=B(X_i,A_r)=0$, and $B(X_i,Z^j)=\delta^j_i$ for all
  $i,j,r$.
 \end{itemize}

 \begin{lemma}\label{lem3.2}
   Adapted frames always exist locally. Moreover, if $\{X_i,A_r,Z^i\}$
   is an adapted local frame, then the dual frame with respect to $B$
   has the form $\{Z^i,A^r,X_i\}$ for certain sections $A^r\in\Ga(\Cal
   A^0M)$. The sections $A_r$ and $A^r$ project to local frames of
   $\Cal A^0M/\Cal A^1M$, which are dual with respect to $B$.
 \end{lemma}
 \begin{proof}
   Consider an open subset $U\subset M$ over which $\Cal G$ is
   trivial.  Then over $U$, each of the bundles $\Cal A^iM$ is trivial
   and hence admits a local frame. We have noted above that $B$
   induces a duality between the bundles $\Cal A^1$ and $\Cal A/\Cal
   A^0$ and a non-degenerate bilinear form on $\Cal A^0/\Cal A^1$.
   Hence we can choose dual local frames $\{Z^i\}$ of $\Cal A^1$ and
   $\{\underline{X}_i\}$ of $\Cal A/\Cal A^0$ as well as an arbitrary
   local frame $\{\underline{A}_r\}$ of $\Cal A^0/\Cal A^1$.  Choosing
   arbitrary representatives $\tilde X_i\in\Ga(\Cal A)$ and $\tilde
   A_r\in\Ga(\Cal A^0)$ for the $\underline{X}_i$ and the
   $\underline{A}_r$, the set $\{\tilde X_i,\tilde A_r,Z^i\}$ forms a
   local frame for $\Cal A$.

   Given such a choice, we put $X_i:=\tilde
   X_i-\sum_j\frac{1}{2}B(\tilde X_i,\tilde X_j)Z^j$. Then
   $B(X_i,X_j)=0$ and $B(X_i,Z^j)=\delta^j_i$ hold by construction.
   Next, def\/ining $A_r:=\tilde A_r-\sum_iB(X_i,\tilde A_r)Z^i$ we also
   get $B(X_i,A_r)=0$, so $\{X_i,A_r,Z^i\}$ is an adapted local frame.

   The form $B$ restricts to zero on $\Cal A^0M\x\Cal A^1M$, so
   $B(Z^i,Z^j)=0$ and $B(Z^i,A_r)=0$ for $i$, $j$, and $r$. This shows
   that the dual frame to $\{X_i,A_r,Z^i\}$ must have the $Z^i$ as its
   f\/irst elements and the $X_i$ as its last elements. Denoting the
   remaining elements by $A^r$, we by def\/inition have $B(A^r,Z^i)=0$
   for all $r$ and $i$. Since $\Cal A^0M$ is the annihilator of $\Cal
   A^1M$ with respect to $B$, this shows that $A^r\in\Ga(\Cal A^0M)$
   for all $r$. The last claim then is obvious.
 \end{proof}

\subsection{A formula for the curved Casimir}\label{3.3}
The key step in our analysis of the curved Casimir operator is that it
admits a simple local formula in terms of an adapted frame. To
formulate this, we need a f\/inal ingredient. The Lie bracket on $\fg$
def\/ines a $P$-equivariant bilinear map $[\ ,\ ]:\fg\x\fg\to\fg$.
Passing to induced bundles, we obtain a skew symmetric bilinear bundle
map $\{\ ,\ \}:\Cal AM\x\Cal AM\to\Cal AM$. By construction, this
makes each f\/iber $\Cal A_xM$ into a f\/iltered Lie algebra isomorphic to
$(\fg,\{\fg^i\})$.
\begin{proposition}\label{prop3.3}
  Let $(p:\Cal G\to M,\om)$ be a parabolic geometry of type $(G,P)$,
  and for a representation $W$ of $P$ consider the natural bundle
  $WM=\Cal G\x_PW$. Let $\Cal C:\Ga(WM)\to\Ga(WM)$ be the curved
  Casimir operator. Let $U\subset M$ be an open subset and
  $\{X_i,A_r,Z^i\}$ be an adapted frame for the adjoint tractor bundle
  $\Cal AM$ over $U$ with dual frame $\{Z^i,A^r,X_i\}$. Then for a
  section $\si\in\Ga(WM)$ we get
\[
\Cal C(\si)|_U=-2\textstyle\sum_i Z^i\bullet
D_{X_i}\si-\sum_i\{Z^i,X_i\}\bullet\si+\sum_rA^r\bullet
A_r\bullet\si.
\]
\end{proposition}
\begin{proof}
  We can use the adapted frame and its dual to evaluate $\Cal C(\si)$
  as
  \begin{equation}
    \label{C1}
    \textstyle\sum_i(D^2\si(X_i,Z^i)+D^2\si(Z^i,X_i))+\sum_r D^2\si(A_r,A^r),
  \end{equation}
  and then use equation \eqref{D^2} from Section~\ref{2.4} to evaluate the
  individual summands. All $Z^i$, $A_r$, and $A^r$ are sections of
  $\Ga(\Cal A^0M)$, so if $s$ is any of these sections then $D_s$ acts
  as $-s\bullet\ $. Naturality of $D$ further implies that for any
  $s\in\Ga(\Cal AM)$, the operator $D_s$ on $\Ga(\Cal AM)$ is
  f\/iltration preserving. Finally, the map $\bullet:\Cal A^0M\x\Cal
  AM\to\Cal AM$ coincides with $\{\ ,\ \}$ and $\{Z^i,X_i\}\in\Ga(\Cal
  A^0)$ for all $i$. Using these facts, the f\/irst sum in \eqref{C1}
  gives
\[
  \textstyle\sum_i\left(-D_{X_i}(Z_i\bullet\si)+(D_{X_i}Z^i)\bullet\si-
    Z_i\bullet D_{X_i}\si-\{Z^i,X_i\}\bullet\si\right),
\]
and by naturality of $D_{X_i}$ the f\/irst two terms add up to
$-Z_i\bullet D_{X_i}\si$. On the other hand, for the last sum in
\eqref{C1}, we obtain
\[
A_r\bullet A^r\bullet\si-\{A_r,A^r\}\bullet\si=A^r\bullet
A_r\bullet\si,
\]
and the result follows.
\end{proof}

\subsection{More basic properties of the curved Casimir}\label{3.4}
One of the main results will be the action of the curved Casimir
operator on bundles induced by irreducible representations. Recall
that $P$ contains the nilpotent normal subgroup $P_+$, and $P/P_+\cong
G_0$. It turns out that $P_+$ acts trivially on any irreducible
representations of $P$, so any such representation comes from an
irreducible representation of the reductive group $G_0$. To formulate
the result, we will need some information on weights.

If $\fg$ is a complex semisimple Lie algebra, then we can choose a
Cartan subalgebra $\frak h\subset\fg$ and positive roots $\Delta^+$ in
such a way that $\frak p\subset\fg$ is a standard parabolic
subalgebra, i.e.~contains $\frak h$ and all positive root spaces.
Then $\frak h\subset\fg_0$ and $\frak h$ is also a Cartan subalgebra
for the reductive Lie algebra $\fg_0$, i.e.~it is the direct sum of
the center $\frak z(\fg_0)$ and a Cartan subalgebra of the semisimple
part of $\fg_0$.  Hence weights for both $\fg$ and $\fg_0$ are linear
functionals on $\frak h$. Our results are easier to formulate in terms
of lowest weights than in terms of highest weights. By $\rho$, we
will denote the lowest form of $\fg$, i.e.~half the sum of all
positive roots, or, equivalently, the sum of all fundamental weights.

Suppose now that $\fg$ is a real semisimple Lie algebra. Then the
grading $\fg=\oplus_{i=-k}^k\fg_i$ induces a grading of the
complexif\/ication $\fg_{\Bbb C}$ and hence $\frak p_{\Bbb
  C}\subset\fg_{\Bbb C}$ is a parabolic subalgebra. Now one can f\/ind a
subalgebra $\frak h\subset\frak g$ such that $ \frak h_{\Bbb
  C}\subset\frak g_{\Bbb C}$ is a Cartan subalgebra for which $\frak
p_{\Bbb C}$ is a standard parabolic. Then we can interpret weights for
$\fg_{\Bbb C}$ as well as for $(\fg_0)_{\Bbb C}$ as linear functionals
on $\frak h$. If~$W$ is a complex irreducible representation of
$\fg_0$, then it extends to $(\fg_0)_{\Bbb C}$, and thus has a lowest
weight. If~$W$ is irreducible without an invariant complex structure,
then the complexif\/ication~$W_{\Bbb C}$ is a complex irreducible
representation of $\fg_0$, and we refer to the lowest weight of
$W_{\Bbb C}$ also as the lowest weight of $W$.

\begin{theorem}\label{thm3.4}
  (1) On any natural vector bundle, the curved Casimir operator is a
      differential operator of order at most one.

(2) Let $W$ be an irreducible representation of $P$. The the curved
    Casimir operator acts by a~scalar on $\Ga(WM)$ for any $(p:\Cal
    G\to M,\om)$. If $-\nu$ is the lowest weight of $W$, then this
    scalar is explicitly given by
\[
\|\nu+\rho\|^2-\|\rho\|^2=\langle\nu,\nu\rangle+2\langle\nu,\rho\rangle,
\]
where the norms and inner products are induced by the Killing form.
\end{theorem}
\begin{proof}
  (1) is obvious from Proposition \ref{prop3.3}, since in the formula for
  $\Cal C(\si)|_U$, the term $D_{X_i}\si$ is f\/irst order in $\si$,
  while all other ingredients are tensorial.

(2) If $W$ is irreducible, then we know that $P_+$ and hence $\frak
p_+$ acts trivially on $W$. Then the formula from Proposition~\ref{prop3.3} reduces to
\[
\Cal C(\si)|_U=-\textstyle\sum_i\{Z^i,X_i\}\bullet\si+\sum_rA^r\bullet
A_r\bullet\si,
\]
so in particular, $\Cal C$ is tensorial and acts only via the values
of $\si$. By construction, this action is induced by the map $W\to W$
given by
\begin{equation}
  \label{Cas}
  w\mapsto -\textstyle\sum_i[Z^i,X_i]\cdot w+\sum_rA^r\cdot
A_r\cdot w,
\end{equation}
where now $\{Z^i\}$ and $\{X_i\}$ are bases of the $\fg_0$-modules
$\frak p_+$ and $\frak
g_-=\fg_{-k}\oplus\dots\oplus\fg_{-1}\cong(\frak p_+)^*$ which are
dual with respect to $B$, while $\{A_r\}$ and $\{A^r\}$ are dual bases
of $\fg_0$ with respect to $B$. In this formula, we only use dual
bases of dual $\fg_0$-modules as well as the brackets $[\ ,\
]:\fg_i\x\fg_{-i}\to\fg_0$ and the action of $\fg_0$ on $W$, which all
are $\fg_0$-equivariant. Thus \eqref{Cas} is a homomorphism of
$\fg_0$-modules, and hence a scalar multiple of the identity by
irreducibility.

To compute this scalar, we may, if necessary, f\/irst complexify $W$ and
then complexify $\fg$. Hence we may without loss of generality assume
that $\fg$ is complex and $W$ is a complex representation of~$\fg$. Now we choose a Cartan subalgebra $\frak h\subset\fg$ such that
$\frak p\subset\fg$ is a standard parabolic subalgebra. Then $\frak h$
and all positive root spaces are contained in $\frak
p=\fg_0\oplus\frak p_+$. For $\al\in\Delta^+$ the root space
$\fg_{\al}$ is thus either contained in $\fg_0$ or in $\frak p_+$ and
in the f\/irst case $\fg_{-\al}\in\fg_0$ while in the second case
$\fg_{-\al}\in\fg_-$. Accordingly, we decompose $\Delta^+$ into the
disjoint union $\Delta^+(\fg_0)\sqcup\Delta^+(\frak p_+)$.

For each root $\al\in\Delta^+$ choose elements $E_\al\in\fg_{\al}$ and
$F_\al\in\fg_{-\al}$ such that $B(E_\al,F_\al)=1$. Then
$[E_\al,F_\al]=H_\al\in H$, the element dual to $\al$ with respect to
$B$. Let $\{H_\ell:\ell=1,\dots,r\}$ be an orthonormal basis of $\frak
h$ with respect to $B$ and put
\begin{gather*}
  \{X_i\}:=\{F_\al:\al\in\Delta^+(\frak p_+)\}\quad
  \{Z^i\}:=\{E_\al:\al\in\Delta^+(\frak p_+)\}\\
\{A_r\}:=\{E_\al,F_\al: \al\in\Delta^+(\frak g_0)\}\cup
\{H_\ell:\ell=1,\dots,r\}.
\end{gather*}
Then these are bases as required, and for the dual basis to $\{A_r\}$
we evidently get
\[
\{A^r\}:=\{F_\al,E_\al: \al\in\Delta^+(\frak g_0)\}\cup
\{H_\ell:\ell=1,\dots,r\}.
\]
Using these bases, we can now evaluate \eqref{Cas} as
\[
  -\sum_{\al\in\Delta^+(\frak p_+)}[E_\al,F_\al]\cdot
  w+\sum_{\al\in\Delta^+(\fg_0)}(E_\al\cdot F_\al+F_\al\cdot
  E_\al)\cdot w+\sum_{\ell=1}^rH_\ell\cdot H_\ell\cdot w.
\]
Each summand in the middle sum can be rewritten as
$(2E_\al\cdot F_\al-[E_\al,F_\al])\cdot w$. Now if $w\in W$ is a
lowest weight vector, then each $F_\al$ acts trivially on $w$, and we
are left with
\[
-\sum_{\al\in\Delta^+}H_\al\cdot w+\sum_{\ell=1}^rH_\ell\cdot
H_\ell\cdot w.
\]
Writing the weight of $w$ as $-\nu$, each summand in the f\/irst sum
produces $\langle\al,\nu\rangle$, and summing up we obtain
$2\langle\rho,\nu\rangle$. Since $\{H_\ell\}$ is an orthonormal
basis, the second sum simply contributes
$\langle-\nu,-\nu\rangle=\|\nu\|^2$.
\end{proof}

\subsection{An application}\label{3.5}
Under weak assumptions, f\/inite dimensional representations of $P$
admit f\/inite composition series by completely reducible
representations.  Correspondingly, any natural bundle has many
irreducible subquotients.  Using the naturality properties of the
curved Casimir operators together with the knowledge of the action on
sections of irreducible bundles, we can construct so-called splitting
operators. These are invariant dif\/ferential operator mapping sections
of some subquotient injectively to sections of the original bundle. We
just describe the general procedure here and do not go into specif\/ic
examples.

Given a representation $W$ of $P$, the composition series is obtained
as follows. For an integer~$N$ to be specif\/ied later, we put
$W^N:=W^{\frak p_+}$, the set of $\frak p_+$-invariant elements of
$W$. Then $W^N\subset W$ is easily seen to be $P$-invariant, and by
construction, $\frak p_+$ acts trivially on $W^N$. Next, one def\/ines
$W^{N-1}:=\{w\in W:Z\cdot w\in W^N\quad\forall\, Z\in\frak p_+\}$. Again,
this is a $P$-invariant subspace of $W$, it contains $W^N$, and
$\frak p_+$ acts trivially on $W^{N-1}/W^N$. Then one continues
inductively and arranges $N$ in such a way that one obtains a
f\/iltration of the form $W=W^0\supset W^1\supset\dots\supset W^N$. If
the center of $\fg_0$ acts diagonalizably on the quotients
$W^i/W^{i+1}$, then all these quotients are completely reducible. By
Theorem~\ref{thm3.4} the curved Casimir then acts on each irreducible
component of this quotient by a scalar that can be computed from the
highest weight.

\begin{theorem}\label{thm3.5}
  Let $W$ be a representation of $P$ endowed with a $P$-invariant
  filtration $W=W^0\supset W^1\supset\dots\supset W^N\supset
  W^{N+1}=\{0\}$ such that each of the subquotients $W^j/W^{j+1}$ is
  completely reducible. Fix a number $i\in\{0,\dots,N\}$ and an
  irreducible component $W'\subset W^i/W^{i+1}$, and let $\mu_0$ be
  the eigenvalue of the curved Casimir operator on the natural bundle
  induced by $W'$. Further, for each $j>i$ let $\mu^j_1,\dots,\mu^j_{n_j}$
  be the different eigenvalues of the curved Casimir operator on the
  bundles induced by the irreducible components of $W^j/W^{j+1}$.

  Let $(p:\Cal G\to M,\om)$ be a parabolic geometry of type $(G,P)$.
  Then the natural operator $L:=\prod_{j=i+1}^N\prod_{k=1}^{n_j}(\Cal
  C-\mu^j_k)$ on $\Ga(WM)$ descends to an operator
  $\Ga(W'M)\to\Ga(W^iM)$. Denoting by $\pi$ the natural tensorial
  projection $\Ga(W^iM)\to\Ga(W^iM/W^{i+1}M)$, we get
  $\pi(L(\si))=(\prod_{j=i+1}^n\prod_{k=1}^{n_j}(\mu^0-\mu^j_k))\si$.
  In particular, if all the $\mu^j_k$ are different from $\mu_0$, then
  $L$ defines a natural splitting operator.
\end{theorem}
\begin{proof}
  For $j=i+1,\dots,N$ put $L_j:=\prod_{k=1}^{n_j}(\Cal C-\mu^j_k)$.
  Then $L=L_{i+1}\o\cdots\o L_N$, and all the operators $L_j$ commute
  with each other. By Proposition \ref{prop3.1}, $\Cal C$ preserves
  sections of each of the subbundles $W^jM$, so in particular, we can
  initially interpret $L$ and the $L_j$ as operators
  $\Ga(W^kM)\to\Ga(W^kM)$ for each $k$.

  Now take $j\in\{i+1,\dots,N\}$ and let $\tilde W\subset W^j/W^{j+1}$
  be an irreducible component. Then we obtain a natural subbundle
  $\tilde WM\subset W^jM/W^{j+1}M$, and we denote by $\pi$ the natural
  projection $W^jM\to W^jM/W^{j+1}M$, as well as the induced tensorial
  operator on sections. For $\si\in\Ga(\tilde WM)$ we can choose
  $\hat\si\in\Ga(W^jM)$ such that $\pi(\hat\si)=\si$. Now
  $\pi^{-1}(\tilde WM)\subset W^jM$ is a natural subbundle, so by
  Proposition \ref{prop3.1}, $\pi(\Cal C(\hat\si))\in\Ga(\tilde WM)$ and
  the image coincides with the action of the curved Casimir of the
  bundle $\tilde WM$ on $\si$. Hence $\pi(\Cal C(\hat\si))=\mu^j_k\si$
  for some $k=1,\dots,n_j$. But this implies that $(\Cal
  C-\mu^j_k)(\hat\si)\in\Ga(W^{j+1}M)$, and since the composition
  def\/ining $L_j$ can be written in arbitrary order, we see that
  $L_j(\hat\si)\in\Ga(W^{j+1}M)$. Since this holds for any irreducible
  component, we see that $L_j$ maps $\Ga(W^jM)$ to $\Ga(W^{j+1}M)$.

  Applied to $j=N$, this says that $L_N$ vanishes on $\Ga(W^NM)$. Thus
  $L_N\o L_{N-1}$ vanishes on $\Ga(W^{N-1}M)$, and inductively, we
  conclude that $L$ vanishes on $\Ga(W^{i+1}M)$. Hence it descends to
  an operator $\Ga(W^iM/W^{i+1}M)\to\Ga(W^iM)$, which we may restrict
  to sections of the natural subbundle $W'M$, so $L$ descends as
  required. In terms of the natural projection
  $\pi:\Ga(W^iM)\to\Ga(W^iM/W^{i+1}M)$ this means that for a section
  $\hat\si\in\Ga(W^iM)$, the value of $L(\hat\si)$ depends only on
  $\pi(\hat\si)$. Supposing that $\pi(\hat\si)\in\Ga(W'M)$,
  Proposition~\ref{prop3.1} and Theorem \ref{thm3.4} imply that we have
  $\pi(L(\hat\si))=(\prod_{j=i+1}^N\prod_{k=1}^{n_j}(\Cal
  C-\mu^j_k))\pi(\hat\si)$, where now $\Cal C$ is the curved Casimir
  on $\Ga(W'M)$. Since this acts by multiplication by $\mu_0$, the
  result follows.
\end{proof}

Let us remark that the number of factors in the def\/inition of $L$ may
sometimes be reduced by passing to the $P$-submodule of $W$ generated
by $W'$.

\section{Tractor bundles and the BGG machinery}\label{4}
In this last section, we analyze the action of the curved Casimir
operator on sections of tractor bundles and on dif\/ferential forms with
values in a tractor bundle. In particular, this will lead to a nice
interpretation of the resolutions of irreducible representations of
$\fg$ by (smooth vectors~in) principal series representations obtained
via BGG sequences.

\subsection{Tractor bundles}\label{4.1}
As before, let $P\subset G$ be a parabolic subgroup in a semisimple
Lie group. Then any representation~$V$ of $G$ can be restricted to
obtain a representation of $P$. The corresponding natural bundles for
parabolic geometries of type $(G,P)$ are called \textit{tractor
  bundles} and play an important role in the theory of parabolic
geometries. The map $G\x V\to G/P\x V$ def\/ined by $(g,v)\mapsto
(gP,g\cdot v)$ induces an isomorphism of homogeneous vector bundles
over the homogeneous model $G/P$, and hence $\Ga(G\x_PV)\cong
C^\infty(G/P,\Bbb R)\otimes V$ as a representation of $G$. However, for a
parabolic geometry $(p:\Cal G\to M,\om)$, the natural bundle $VM=\Cal
G\x_PV$ is in general not trivial.

Via the trivialization constructed above, the bundle $G\x_PV$ inherits
a canonical $G$-invariant linear connection. Such a natural
connection on a tractor bundle also exists for arbitrary parabolic
geometries, which is the main reason for the importance of these
bundles. The general theory of these tractor connections is developed
in \cite{TAMS}. For our purposes, the simplest description of the
tractor connection is the following. Since $V$ is a representation of
$G$, the inf\/initesimal action of $\frak p$ extends to $\frak g$. Hence
the bundle map $\bullet$ introduced in Section~\ref{2.3} extends to a bundle
map $\Cal AM\x VM\to VM$, which has the properties of a Lie algebra
action.

Now consider the operator $\Ga(\Cal AM)\x\Ga(VM)\to\Ga(VM)$ def\/ined by
$(s,\si)\mapsto D_s\si+s\bullet\si$. This is linear over
$C^\infty(M,\Bbb R)$ in the $\Cal AM$-slot and satisf\/ies a Leibniz
rule in the $VM$-slot. From Section~\ref{2.3} we know that for $s\in\Ga(\Cal
A^0M)$, we get $D_s\si=-s\bullet\si$, so our bundle map descends to
$\Ga(\Cal AM/\Cal A^0M)$ in the f\/irst argument. By construction, $\Cal
AM/\Cal A^0M=\Cal G\x_P\frak g/\frak p\cong TM$. Hence our bundle map
descends to a bundle map $\frak X(M)\x\Ga(VM)\to\Ga(VM)$, which is
linear over smooth functions in the vector f\/ield and satisf\/ies a
Leibniz rule in the other slot. Thus it def\/ines a natural linear
connection $\nabla^V$ on $VM$.

We have noted already that the natural bundle corresponding to $\frak
g/\frak p$ is the tangent bundle. Moreover, from Section~\ref{3.2} we know
that $(\frak g/\frak p)^*\cong\frak p_+$ as a $P$-module, so the
natural bundle corresponding to $\frak p_+$ is the cotangent bundle.
In particular, for a parabolic geometry on $M$, the cotangent bundle
$T^*M$ naturally is a bundle of nilpotent f\/iltered Lie algebras
modeled on $\frak p_+$. Moreover, the cotangent bundle $T^*M$ can be
naturally identif\/ied with the subbundle $\Cal A^1M$ of the adjoint
tractor bundle. Given a representation $W$ of $P$, the bundles of
dif\/ferential forms with values in the natural bundle $WM$ correspond
to the representations $\La^k\frak p_+\otimes W$.

Now this gives rise to another natural operation. Via the
inf\/initesimal action, $\frak p_+$ acts on any representation $W$ of
$P$. The corresponding bundle map is simply the restriction
$\bullet:\Cal A^1M\x WM\to WM$ of the bundle map introduced in Section~\ref{2.3}. Otherwise put, $WM$ is a bundle of modules over the bundle
$T^*M$ of Lie algebras. Now the Lie algebra homology dif\/ferential
def\/ines a natural map $\La^k\frak p_+\otimes W\to\La^{k-1}\frak
p_+\otimes W$. In the literature on parabolic geometries, this
operation is mainly considered in the case of a representation of
$\fg$ and it is called the \textit{Kostant codifferential} and denoted
by $\partial^*$, since it was introduced in \cite{Kostant} as the
adjoint of a Lie algebra cohomology dif\/ferential. Explicitly, the
corresponding bundle map $\partial^*:\La^kT^*M\otimes
WM\to\La^{k-1}T^*M\otimes WM$ is on decomposable elements given by
  \begin{gather*}
    \partial^*(\ph_1\wedge\dots\wedge\ph_k\otimes\si):=\sum_{i=1}^k(-1)^i\ph_1
\wedge\dots\wedge\widehat{\ph_i}\wedge\dots\wedge\ph_k\otimes\ph_i\bullet\si\\
\phantom{\partial^*(\ph_1\wedge\dots\wedge\ph_k\otimes\si):=}{} +\sum_{i<j}(-1)^{i+j}\{\ph_i,\ph_j\}\wedge\ph_1\wedge\dots\wedge
\widehat{\ph_i}\wedge\dots\wedge\widehat{\ph_j}\wedge\dots\wedge
\ph_k\otimes\si.
  \end{gather*}

Now we have all ingredients at hand to compute the action of the
curved Casimir on a tractor bundle.
\begin{proposition}\label{prop4.1}
  Let $V$ be a finite dimensional irreducible representation of $G$,
  and let $c_0$ be the scalar by which the Casimir operator acts on
  $V$. Let $(p:\Cal G\to M,\om)$ be a parabolic geometry of type
  $(G,P)$, let $VM\to M$ be the tractor bundle corresponding to $V$,
  and let $\nabla^V$ be the canonical tractor connection on this
  bundle.

  Then the curved Casimir operator $\Cal C:\Ga(VM)\to\Ga(VM)$ is given
  by $\Cal C(\si)=c_0+2\partial^*\nabla^V\si$.
\end{proposition}
\begin{proof}
  It suf\/f\/ices to verify this locally over an open subset $U\subset M$
  over which there is an adapted frame $\{X_i,A_r,Z^i\}$ for $\Cal
  AM$, and use the formula from Proposition \ref{prop3.3}. Since we are
  working on a tractor bundle, we can expand
\[
-\{Z^i,X_i\}\bullet\si=-Z^i\bullet X_i\bullet\si+X_i\bullet Z^i\bullet\si.
\]
Using the relation between the tractor connection and the fundamental
derivative, we conclude that $\Cal C(\si)|_U$ is given by
\[
-2\sum_iZ^i\bullet\nabla_{\Pi(X_i)}^V\si+\sum_i(X_i\bullet
Z^i\bullet\si+Z^i\bullet X_i\bullet\si)+\sum_rA^r\bullet
A_r\bullet\si,
\]
where $\Pi:\Cal AM\to TM$ is the natural projection. By
construction, $\{\Pi(X_i)\}$ is a local frame for~$TM$ over $U$ and
$\{Z^i\}$ is the dual frame of $T^*M$, so the f\/irst sum in this
formula gives~$2\partial^*\nabla^V\si$. The rest of the formula is
tensorial, with the elements only acting on the values of $\si$. But
on the values, we simply obtain the iterated actions of the elements
of dual bases of $\fg$, so this produces the action of the Casimir on
$V$.
\end{proof}

\subsection{Tractor bundle valued forms}\label{4.2}
The canonical connection $\nabla^V$ naturally extends to the
\textit{covariant exterior derivative} $d^{\nabla^V}$ on tractor bundle
valued dif\/ferential forms. The def\/inition of
$d^{\nabla^V}:\Om^k(M,VM)\to\Om^{k+1}(M,VM)$ is obtained by replacing
in the standard formula for the exterior derivative the actions of
vector f\/ields by the covariant derivative with respect to $\nabla^V$.
It is well known that $d^{\nabla^V}\o d^{\nabla^V}$ is given by the
action of the curvature $R^V$ of $\nabla^V$. For a tractor connection,
this curvature is induced by the curvature of the canonical Cartan
connection, see Proposition 2.9 of \cite{TAMS}. In particular, for the
homogeneous model, $d^{\nabla^V}\o d^{\nabla^V}=0$, and
$(\Om^*(G/P,G\x_PV),d^{\nabla^V})$ is a f\/ine resolution of the
constant sheaf $V$ on $G/P$. This is called the \textit{twisted
  de-Rham resolution}.

On the other hand, there is a natural resolution of the f\/inite
dimensional representation $V$ by generalized Verma modules. This is
the \textit{(generalized) Bernstein--Gelfand--Gelfand (BGG)
  resolution}, constructed in \cite{BGG} for the Borel subgroup and in
\cite{Lepowsky} in general. The construction of these resolutions is
purely algebraic, using existence results for homomorphisms between
generalized Verma modules, and combinatorial properties of Weyl
groups. Now there is a duality between homomorphisms of generalized
Verma modules and invariant dif\/ferential operators between sections of
homogeneous bundles, see appendix A of the preprint version of
\cite{CSS-BGG} for an exposition. Via this duality, the sequence of
homomorphisms forming the BGG resolution corresponds to a~sequence of
invariant dif\/ferential operators on sections of homogeneous bundles
over $G/P$.

As mentioned before, constructing invariant dif\/ferential operators for
parabolic geometries is an important and dif\/f\/icult problem. It is
known that the BGG resolutions contain a majority of the existing
homomorphisms between generalized Verma modules, so trying to obtain
the corresponding operators for curved geometries is a natural idea.
This problem was solved in~\cite{CSS-BGG} and the construction was
signif\/icantly improved in~\cite{David-Tammo}. The idea is to construct
the BGG operators from the twisted de-Rham resolution in a way which
extends to the curved case. There the composition of two covariant
exterior derivatives is non-zero, so the twisted de-Rham sequence is
not a complex any more. Hence one only obtains a sequence of invariant
dif\/ferential operators, called a BGG-sequence. It should be mentioned
however, that under certain semi-f\/latness assumptions, curved BGG
sequences may contain very interesting subcomplexes, see~\cite{subcomplexes}.

The relation between the two sequences is provided by Lie algebra
homology. Let $(p:\Cal G\to M,\om)$ be a parabolic geometry of type
$(G,P)$, let $V$ be a representation of $G$ and consider the operators
$\partial^*$ acting on $\La^kT^*M\otimes VM$ as def\/ined in Section~\ref{4.1}.
It is well known that $\partial^*\o\partial^*=0$, so
$\im(\partial^*)\subset\ker(\partial^*)\subset\La^kT^*M\otimes VM$ are
natural subbundles. The subquotient $\ker(\partial^*)/\im(\partial^*)$
is denoted by $H_k(T^*M,VM)$, since it the just the f\/iber-wise Lie
algebra homology. By construction it is the associated bundle $\Cal
G\x_PH_k(\frak p_+,V)$. A crucial fact (see \cite{CSS-BGG}) is that
the natural action of~$\frak p_+$ on $\La^k\frak p_+\otimes V$ maps
$\ker(\partial^*)$ to $\im(\partial^*)$. Hence $\frak p_+$ acts
trivially on the homology $H_k(\frak p_+,V)$, so this is a completely
reducible representation.

All this also works for a general representation $W$, but then it is
not clear what the homology groups look like. However, for a
representation $V$ of $\frak g$, the $\fg_0$-representations
$H_k(\frak p_+,V)$ are computed by Kostant's version of the
Bott--Borel--Weil theorem proved in \cite{Kostant}. Kostant actually
computes the Lie algebra cohomology groups $H^k(\frak p_+,V)$, but
they are dual to the homology groups. The description is algorithmic,
so the computation can be easily carried out explicitly. The result is
that these are exactly the duals of the representations showing up in
the Verma-module version of the BGG resolution, so they are exactly
the representations we need.

By construction, for each $k$ and the subbundle
$\ker(\partial^*)\subset \La^kT^*M\otimes VM$, we have a natural
bundle projection $\pi_H:\ker(\partial^*)\to H_k(M,VM)$. This induces
a tensorial operator on sections, which we denote by the same symbol.
The crucial step in the construction of BGG sequences is to construct
invariant splitting operators
$L:\Ga(H_k(M,VM))\to\Ga(\ker(\partial^*))$ such that $\pi_H\o L=\id$.
This is done by looking at
$\square^R:=d^{\nabla^V}\o\partial^*+\partial^*\o d^{\nabla^V}$, which
is an invariant f\/irst order operator def\/ined on $\Om^k(M,VM)$. By
construction, $\square^R$ commutes with $\partial^*$ and hence
preserves sections of the natural subbundles $\ker(\partial^*)$ and
$\im(\partial^*)$. Further, it turns out that, similarly as discussed
for the adjoint tractor bundle in Section~\ref{3.2}, any tractor bundle is
canonically f\/iltered. The operator $\square^R$ is f\/iltration
preserving, i.e.~maps sections of any f\/iltration component to sections
of the same component.

Hence one can look at the lowest homogeneous component, or otherwise
put, the induced operator on the associated graded vector bundle. This
turns out to be tensorial and induced by the Kostant--Laplacian
$\square$, which plays a crucial role in Kostant's proof in
\cite{Kostant}. In particular, the action of $\square$ on sections of
the subbundle $\im(\partial^*)$ is immediately seen to be invertible.
From this, one concludes that $\square^R$ acts invertibly on
$\Ga(\im(\partial^*))$, and the inverse $Q$ is a dif\/ferential operator,
which is invariant by construction. Then for $\si\in\Ga(H_k(M,VM))$
one can choose any $\ph\in\Ga(\ker(\partial^*))$ such that
$\pi_H(\ph)=\si$ and def\/ine $L(\si):=\ph-Q(\square^R(\ph))$. It turns
out that $d^{\nabla^V}\o L$ has values in $\ker(\partial^*)$, so one
may def\/ine the BGG operators as $D^V:=\pi_H\o d^{\nabla^V}\o L$.
Finally, for $\ph\in\Ga(\ker(\partial^*))$, vanishing of
$\partial^*d^{\nabla^V}\ph$ already implies vanishing of
$\square^R(\ph)$ and hence $\ph=L(\pi_H(\ph))$, which is a very useful
characterization of the splitting operator.

\subsection{The relation to the curved Casimir}\label{4.3}
We next discuss the relation of the operator $\square^R=\partial^*\o
d^{\nabla^V}+d^{\nabla^V}\o\partial^*$ to the curved Casimir operator.
We will get a complete description in the locally f\/lat case and a
large amount of information in general. The main tool for these
results is the following lemma which holds for general geometries.
\begin{lemma}\label{lem4.3}
  Let $(p:\Cal G\to M,\om)$ be a parabolic geometry of type $(G,P)$,
  $V$ an irreducible representation of $G$, and $VM$ the corresponding
  tractor bundle. Let $c_0\in\Bbb R$ be the eigenvalue of the Casimir
  operator of $\frak g$ on the irreducible representation $V$.

  Then the action of the curved Casimir operator $\Cal C$ on
  $\Om^k(M,VM)$ is compatible with the natural filtration on
  $\La^kT^*M\otimes VM$ discussed in Section~\ref{4.2} and the induced action
  on the sections of the associated graded vector bundle is given by
  $2\square+c_0$, where $\square$ denotes the bundle map induced by the
  Kostant--Laplacian.
\end{lemma}
\begin{proof}
  The f\/iltration of $\La^kT^*M\otimes VM$ is formed by natural
  subbundles, so the fact that sections of these subbundles are
  preserved by $\Cal C$ follows immediately from naturality as proved
  in Proposition~\ref{prop3.1}. The associated graded vector bundle is
  induced by a representation with trivial action of $\frak p_+$ and
  diagonalizable action of $\frak z(\frak g_0)$, so it splits into a
  direct sum of bundles induced by irreducible representations of
  $P$. Again by naturality, the operator induced on the associated
  graded by $\Cal C$ coincides with the curved Casimir on the bundle
  in question.

  Now the proof of Theorem \ref{thm3.4} shows that the latter curved
  Casimir acts by a scalar on each lowest weight vector, and hence on
  each isotypical component. Using the conventions of Section~\ref{3.4} on
  weights, this scalar is explicitly given by
  $\|\nu+\rho^2\|-\|\rho\|^2$ on the bundle corresponding to the
  isotypical component of lowest weight $-\nu$.

  Now it is a crucial part of Kostant's proof in \cite{Kostant}, that
  the Kostant--Laplacian acts by a scalar on each isotypical component.
  Carrying out the dualizations needed to translate from Kostant's
  setting (cohomology of $\frak p_+$) to the one we need here
  (homology of $\frak p_+$ with coef\/f\/icients in the dual
  representation), this scalar is given by
  $\frac{1}{2}(\|\nu+\rho\|^2-\|\la+\rho\|^2)$ on the isotypical
  component of weight $-\nu$, where $-\la$ is the lowest weight of
  $V$. But then the claim follows immediately since if $-\la$ is the
  lowest weight of $V$, then $c_0=\|\la+\rho^2\|-\|\rho\|^2$.
\end{proof}

On the one hand, this result immediately implies that the operator
$\Cal C-c_0$ on $\Om^k(M,VM)$ is f\/iltration preserving and induces
$2\square$ on the associated graded. Hence one can construct curved BGG
sequences replacing the operator $\square^R$ by $\Cal C-c_0$. There is
not much dif\/ference between these two operators, however, because we
have the following result.
\begin{proposition}\label{prop4.3}
  On any parabolic geometry $(p:\Cal G\to M,\om)$ of type $(G,P)$ and
  any tractor bundle~$VM$, the operator $2\square^R-\Cal C+c_0$ on
  $\Om^k(M,VM)$ is tensorial and maps each component of the natural
  filtration of $\La^kT^*M\otimes VM$ to the next smaller filtration
  component.
\end{proposition}
\begin{proof}
  Since both $2\square^R$ and $\Cal C-c_0$ are f\/iltration preserving,
  so is their dif\/ference. But by Lemma~\ref{lem4.3}, this dif\/ference induces
  the zero map on the associated graded bundle, and hence moves each
  f\/iltration component to the next smaller one. So it remains to prove
  that the dif\/ference is tensorial, i.e.~linear over $C^\infty(M,\Bbb
  R)$.

Now $\partial^*$ is tensorial, and for $d^{\nabla^V}$ we have the
Leibniz rule $d^{\nabla^V}(f\ph)=df\wedge\ph+fd^{\nabla^V}\ph$. This
implies that
\[
\square^R(f\ph)=f\square^R(\ph)+\partial^*(df\wedge\ph)+
df\wedge\partial^*\ph= f\square^R(\ph)-df\bullet\ph,
\]
where for the last step we used formula (1.2) of \cite{CSS-BGG}.

On the other hand, the fundamental derivative satisf\/ies a Leibniz
rule by Proposition~3.1 of~\cite{TAMS} and naturality. Moreover, on
smooth functions, the fundamental derivative coincides with the
exterior derivative and the action map $\bullet$ is trivial. Using
these facts, the formula for $\Cal C$ in terms of adapted local frames
from Proposition~\ref{prop3.3} immediately implies that
\[
\Cal C(f\ph)=f\Cal C(\ph)-2\textstyle\sum_idf(X_i)Z^i\bullet\ph=f\Cal
C(\ph)-2df\bullet\ph.
\]
Together with the above, this shows that $2\square^R-\Cal C+c_0$ is
tensorial, which completes the proof.
\end{proof}

Using this result, we can completely clarify the situation on the
homogeneous model, and hence for locally f\/lat geometries (which are
locally isomorphic to the homogeneous model).
\begin{corollary}\label{cor4.3}
  Suppose that $V$ is an irreducible representation of $\fg$ and let
  $c_0$ be the eigenvalue of the Casimir operator on $V$. Then on the
  homogeneous model and hence on locally flat geometries, the curved
  Casimir operator on differential forms with values in the tractor
  bundle induced by~$V$ is given by $\Cal C=2\square^R+c_0$.
\end{corollary}
\begin{proof}
  From Proposition~\ref{prop4.3} above we know that for each $k$, the operator
  $2\square^R+c_0-\Cal C$ is a~tensorial invariant operator on
  $G\x_P(\La^k\frak p_+\otimes V)$, so this is induced by a
  $P$-equivariant linear map on $\La^k\frak p_+\otimes V$. In
  addition, the proposition tells us that this linear map has to map
  each of the components of the natural $P$-invariant f\/iltration of
  $\La^k\frak p_+\otimes V$ into the next smaller f\/iltration
  component. But it is well known that there is an element $A\in\frak
  g_0$ (called the grading element) which acts diagonalizably on
  $\La^k\frak p_+\otimes V$ in such a way that the eigenvalues are of
  the form $a_0,a_0+1,\dots,a_0+N$ and the f\/iltration components are
  just the sum of the eigenspaces for the eigenvalues
  $a_0+i,\dots,a_0+N$ for some $i$. But then equivariancy of our map
  for the action of $A$ implies that it is identically zero.
\end{proof}

\subsection{An interpretation of the BGG-resolution}\label{4.4}
Let $V$ be an irreducible representation of $\fg$ and consider the
twisted de-Rham resolution \linebreak $(\Om^*(G/P, G\x_PV),d^{\nabla^V})$
discussed in Section~\ref{4.2} on the homogeneous model. In terms of
representation theory, this is a resolution of the locally constant
sheaf $V$ on $G/P$ by induced representations and $G$-equivariant
(dif\/ferential) operators. On each of the spaces $\Om^k(G/P,G\x_PV)$ we
have the Casimir operator $\Cal C=\Cal C_k$, and naturality of the
Casimir implies that $d^{\nabla^V}\o\Cal C_k=\Cal C_{k+1}\o
d^{\nabla^V}$. In particular, the operators $d^{\nabla^V}$ map
eigenspaces of the Casimir operator $\Cal C$ to eigenspaces for the
same eigenvalue. The relevant eigenvalue here of course is the
eigenvalue~$c_0$ on the representation $V$, which is isomorphic to the
space of sections of $G\x_PV$ which are parallel for $\nabla^V$. Using
known results on BGG sequences, we can determine the eigenspace of the
Casimir for the eigenvalue $c_0$.

\begin{proposition}\label{prop4.4}
  The projection $\pi_H:\Om^k(G/P,G\x_PV)\to\Ga(G\x_PH_k(\frak
  p_+,V))$ restricts to a linear isomorphism on the eigenspace of the
  Casimir operator $\Cal C$ for the eigenvalue $c_0$. In particular,
  this eigenspace is a finite direct sum of (the spaces of all smooth
  vectors in) principal series representations.
\end{proposition}
\begin{proof}
  From Corollary~\ref{cor4.3} we know that the $c_0$-eigenspace of $\Cal
  C$ coincides with the kernel of $\square^R=\partial^*\o
  d^{\nabla^V}+d^{\nabla^V}\o\partial^*$. Now let $L: \Ga(G\x_PH_k(\frak
  p_+,V))\to\Om^k(G/P,G\x_PV)$ be the BGG splitting operator as
  described in Section~\ref{4.2}. Then $\pi_H\o L=\id$, so $\pi_H$ is
  surjective. We have noted in Section~\ref{4.2} that $\partial^*\o L=0$ and
  $\partial^*\o d^{\nabla^V}\o L=0$, and hence $L$ has values in
  $\ker(\square^R)$.

  Conversely, assume that $\square^R(\ph)=0$. By construction,
  $\square^R$ commutes with $\partial^*$, so this implies
  $\square^R(\partial^*\ph)=0$. But we have noted in Section~\ref{4.2} that
  $\square^R$ acts invertibly on sections of $\im(\partial^*)$. Hence
  $\partial^*\ph=0$ and we can form $\pi_H(\ph)$. But since
  $\partial^*\ph=0$, we get
  $0=\square^R(\ph)=\partial^*d^{\nabla^V}\ph$, and in Section~\ref{4.2} we
  have noted that these facts imply that $\ph=L(\pi_H(\ph))$. Hence we
  conclude that~$\pi_H$ and~$L$ restrict to inverse isomorphisms
  between $\ker(\square^R)$ and $\Ga(G\x_PH_k(\frak p_+,V))$.
\end{proof}

In Theorem~2.6 of~\cite{CSS-BGG} it is shown that on the homogeneous
model, any BGG sequence is a~complex, which def\/ines a f\/ine resolution
of the f\/inite dimensional representation~$V$. Here we obtain an
alternative proof of this fact as well as a nice interpretation in
terms of representations theory.

\begin{corollary}\label{cor4.4}
  Restricting the twisted de-Rham resolution
  $(\Om^*(G/P,G\x_PV),d^{\nabla^V})$ to the eigen\-spaces of the
  Casimir operator for the eigenvalue $c_0$ induces a finite
  resolution of the representation $V$ by direct sums of (spaces of
  smooth vectors in) principal series representations. The inducing
  $P$-representations for these principal series representations are
  the duals of the representations inducing the generalized Verma
  modules in the BGG resolution.
\end{corollary}
\begin{proof}
  It is an elementary exercise to verify that the restriction of the
  twisted de-Rham resolution to eigenspaces for one f\/ixed eigenvalue
  of the Casimir is again a resolution of the eigenspace in the kernel
  of the f\/irst operator for the given eigenvalue. Choosing the
  eigenvalue to be $c_0$, this eigenspace is (isomorphic to) $V$. By
  Proposition~\ref{prop4.4}, the $c_0$-eigenspace in each $\Om^k(G/P,G\x_PV)$
  is a direct sum as required. The inducing representations are
  described by Kostant's version of the Bott--Borel--Weil theorem,
  which proves the last claim.
\end{proof}

\subsection*{Acknowledgments} The idea to study the Casimir operator
on tractor bundle valued forms grew out of questions by M.~Cowling and
by P.~Julg, who conjectured Corollary~\ref{cor4.3} for the case of the
trivial representation. We are very grateful to them for drawing our
attention to this problem. Our thanks also go to the anonymous
referees for helpful suggestions and corrections. Most of the work was
done during meetings of authors at the Erwin Schr\"odinger Institute
for Mathematical Physics in Vienna. First author supported by project
P19500--N13 of the Fonds zur F\"orderung der wissenschaftlichen
Forschung (FWF). The second author thanks the grant GA\v CR
Nr.~201/05/2117 and the institutional grant MSM 0021620839 for their
support.

\pdfbookmark[1]{References}{ref}
\LastPageEnding

\end{document}